\documentclass[12pt]{article}
\usepackage{amssymb,amsmath, amsthm,latexsym, verbatim, amscd}
\usepackage{color}

\usepackage{geometry}
\usepackage[all]{xy}
\usepackage{graphicx}
\usepackage{xr}

\newtheorem{theorem}{Theorem}[section]
\newtheorem{prop}[theorem]{Proposition}
\newtheorem{lemma}[theorem]{Lemma}
\newtheorem{cor}[theorem]{Corollary}
\theoremstyle{definition}
\newtheorem{definition}[theorem]{Definition}
\newtheorem{example}[theorem]{Example}
\newtheorem{notation}[theorem]{Notation}
\theoremstyle{remark}
\newtheorem{remark}[theorem]{Remark}

\newcommand{\btheorem}{\begin{theorem}}
\newcommand{\etheorem}{\end{theorem}}
\newcommand{\bprop}{\begin{prop}}
\newcommand{\eprop}{\end{prop}}
\newcommand{\blemma}{\begin{lemma}}
\newcommand{\elemma}{\end{lemma}}
\newcommand{\bcor}{\begin{cor}}
\newcommand{\ecor}{\end{cor}}

\numberwithin{figure}{section}

\newcommand{\Exterior}{\mathchoice{{\textstyle\bigwedge}}%
    {{\bigwedge}}%
    {{\textstyle\wedge}}%
    {{\scriptstyle\wedge}}} 

\newcommand{\bg}{{\mathfrak g }}

\newcommand{\bq}{{\mathfrak q }}

\newcommand{\bC}{{\mathbb C}}
\newcommand{\bZ}{{\mathbb Z}}

\title{
Distinguished representations of 
\\
 $SO(n+1,1) \times SO(n,1)$, periods and branching laws}
\author
{
Toshiyuki Kobayashi
\vspace{-10mm}
\footnote
{
Graduate School of Mathematical Sciences,
 The University of Tokyo
 and Kavli IPMU (WPI),
 3-8-1 Komaba, Tokyo, 153-8914 Japan
\newline
\textit{E-mail address}:
\texttt{toshi@ms.u-tokyo.ac.jp}
}
\,\,
 and 
 Birgit Speh
\vspace{-10mm}
\footnote
{Department of Mathematics, Cornell University, Ithaca, NY 14853-4201, USA
\newline
\textit{E-mail address}:
\texttt{bes12@cornell.edu}}
}

\begin{document}

\maketitle

\noindent
\begin{abstract}
Given irreducible representations $\Pi$ and $\pi$
 of the rank one special orthogonal groups $G=SO(n+1,1)$ and $G'=SO(n,1)$
 with nonsingular integral infinitesimal character, 
 we state in terms of $\theta$-stable parameter
 necessary and sufficient conditions so that 
\[   
\operatorname{Hom}_{G'}(\Pi|_{G'}, \pi )\not = \{0\}.    
\]
In the special case that both  $\Pi$ and $\pi$ are tempered,
 this implies the Gross--Prasad conjectures for tempered representations
 of $S O(n+1,1) \times S O(n,1)$ which are nontrivial on the center.

We apply these results to construct nonzero periods and distinguished representations. 
If both $\Pi$ and $ \pi$ have the trivial infinitesimal character $\rho$
 then we  use a theorem
 that the periods are nonzero on the minimal $K$-type
 to obtain a nontrivial bilinear form on the $(\bg,K)$-cohomology 
 of the representations.  
\end{abstract}

\textbf{Mathematic Subject Classification (2010):} 
Primary 22E30; 
Secondary 11F70, 
22E45, 
22E46, 
53A30, 
58J70

\section{Introduction}

In the book \cite{Weyl97}
{\it The Classical Groups: Their Invariants and Representations}
 published in 1939, 
 Hermann Weyl
 discusses the restriction of irreducible finite-dimensional representations of the orthogonal group $O(n)$ to an orthogonal subgroup $O(n-1)$.  
The restriction splits into a direct sum of irreducible representations and can be roughly described by the following {\it branching law.}
Suppose $n \ge 2$. 
Denote by $[x]$ the greatest integer
 that does not exceed $x$.  
To each irreducible representation of $O(n)$ is assigned a highest weight 
\[   
   \mu =
   (\mu_1,\mu_2, \dots ,\mu_{   [\frac n 2]} ) 
\] 
satisfying
\[   
   \mu_1 \ge \mu_2 \ge \dots  \ge  \mu_{   [\frac n 2]}   
\]
and a character $\varepsilon$ of $O(n)/SO(n)$.  
We denote by this representation $ F^{O(n)} ( \mu)_\varepsilon $ where $\varepsilon \in \{ +,-\}$
 ($\varepsilon$ is not unique when $n$ is even and $\mu_{\frac n 2} \ne 0$).  
H. Weyl obtained the \lq\lq{branching law}\rq\rq\ as

\begin{equation}
\label{eqn:branchWeyl} 
   F^{O(n)} (\mu_1,  \dots, \mu_{[\frac n 2]})_\varepsilon  |_{O(n-1) }
 = \bigoplus F^{O(n-1)} (\nu_1,\dots, \nu_{[\frac{n-1}{2}]}   ) _\varepsilon
\end{equation}
where the summation is taken over 
$(\nu_1,\dots, \nu_{[\frac{n-1}{2}]}) \in \bZ^{[\frac{n-1}{2}]}$
 subject to 
\begin{alignat*}{2}
&\mu_{1} \ge \nu_{1} \ge \mu_{2} \ge \cdots \ge \nu_{\frac{n-1}{2}} \ge 0
\quad
&&\text{for $n$ odd}, 
\\
&\mu_{1} \ge \nu_{1} \ge \mu_{2} \ge \cdots \ge \nu_{\frac{n-2}{2}}
\ge \mu_{\frac {n} 2}\ge 0
\quad
&&\text{for $n$ even}.   
\end{alignat*}

\bigskip
In this article
 we present similar branching laws
 for the restriction
 of irreducible {\it{infinite-dimensional}} representations
 of $SO(n+1,1)$ to the subgroup $SO(n,1)$. 
Since the restriction of an infinite-dimensional representation $\Pi$ of $SO(n+1,1)$ to $SO(n,1)$
 is not a direct sum of irreducible representations \cite{xkInvent98}, 
 we consider as in \cite{sbon, sbonvec} the representations $\Pi$ and $\pi$ of $SO(n+1,1)$,
 respectively $SO(n,1)$, 
 realized as smooth representations of moderate growth \cite[Chap.~11]{W}
 and define the {\it multiplicity} by
\[
   m(\Pi,\pi) := \mbox{dim}_{\mathbb{C}} \operatorname{Hom}_{SO(n,1)}  (\Pi|_{SO(n,1)}, \pi).
\]
The multiplicity is either 0 or 1 \cite{SunZhu}.

We consider in the article only  representations of the special orthogonal group which  have the same infinitesimal character 
 as an irreducible  finite-dimensional representation $F$ of $S O(n+1)$.  
To simplify the notation and presentation we assume in this article
 that $F$ is \lq\lq{self-dual}\rq\rq, 
{\it{i.e.}}, 
 assume that the highest weight 
 $\mu=(\mu_1, \cdots, \mu_{m+1})$ satisfies $\mu_{m+1} = 0$ when $n=2m+1$.  
We do not impose any assumption
 on $F$ when $n$ is even.  
See Assumption A in Section \ref{subsec:II.1.1}.  
For the general case,
 see \cite{sbonvec} and \cite{branching}.

For every irreducible representation $\pi$ of $S O(n,1)$
 we define in Section \ref{sec:II.3} a {\it{height}}  
\[ 
   h_\pi \in \{0,\dots m \}  \quad \mbox{if }
  \,\,\text{$n=2m$ or $2m+1$}  
\]
and in Section \ref{subsec:II.4} a {\it{signature}} $\delta \in \{+,-\}$. 
The signature is unique except for discrete series representations. 
If $\pi$  has the same infinitesimal character as $F(\mu)$ 
 we say
 that $(\mu, h_{\pi}, \delta)$ are
 the {\it{enhanced $\theta$-stable parameters}} of $\pi$. 
The representations with {\emph {enhanced $\theta$-stable parameters}}
 $((0,\dots,0), i, \delta)$ are representations
 with trivial infinitesimal character $\rho$
 and are denoted by $\Pi_{i,\delta}$.
See \cite[Chap.~2, Sect.~4]{sbonvec} for a description of the representations
 of $O(n+1,1)$
 with trivial infinitesimal character $\rho$
 and \cite[Chap.~14, Sect.~9]{sbonvec} 
 for their enhanced $\theta$-parameters.

\begin{remark}
\label{rem:I.1}
In \cite{sbonvec}
we have treated mainly the full group $O(n+1,1)$
 rather than the special orthogonal group $S O(n+1,1)$, 
 and stated results for $S O(n+1,1)$
 in Chapter 11 with \lq\lq{bar}\rq\rq\
 for the corresponding objects.  
The relation between branching laws 
 for $O(n+1,1) \downarrow O(n,1)$ 
 and $S O(n+1,1) \downarrow S O(n,1)$ is discussed 
 in \cite[Chap. 15 (Appendix II)]{sbonvec}.  
In this article, 
 we treat mainly the special orthogonal group 
 $S O(n+1,1)$,
 and use different convention 
 in the point that we omit the \lq\lq{bar}\rq\rq\ 
 for representations of $S O(n+1,1)$.  
\end{remark}

\medskip
All results in this article  are based
 on the following {\bf  branching theorem}:

\begin{theorem}[branching law]
\label{branch}
Let  $\Pi$ and $\pi$  be irreducible representations of $SO(n+1,1)$ respectively $SO(n,1)$ with enhanced $\theta$-stable parameters $(  \mu, h_{\Pi,} \varepsilon)$, respectively $(\nu, h_{\pi}, \delta)$.
\begin{enumerate}
\item[{\rm{(1)}}]
Suppose  that $n=2m$.  
Then 
\[
\operatorname{dim}_{\mathbb{C}}\operatorname{Hom}_{SO(n,1)}(\Pi|_{SO(n,1)}, \pi)=1
\]
 if and only if 
the enhanced $\theta$-stable parameters of the representations $\Pi$ and $\pi $ satisfy 
\begin{enumerate}
\item[{\rm{(a)}}]  $\varepsilon=\delta$, 
\item[{\rm{(b)}}]
$h_\pi \in \{ h_\Pi, h_\Pi -1\}$ when $h_\pi < m$ 
and 
$h_\pi =h_\Pi$ when $h_\pi =m $, 
\item[{\rm{(c)}}]  
$
 \mu_0 \geq \nu_0 \geq \mu_1 \geq \dots \geq \mu_m  \geq \nu_m \geq 0$.  
\end{enumerate}

\item[{\rm{(2)}}] Suppose  that $n=2m-1$. 
Then 
\[
   \operatorname{dim}_{\mathbb{C}}
   \operatorname{Hom}_{SO(n,1)}(\Pi|_{SO(n,1)}, \pi) =1
\]
 if and only if 
the enhanced $\theta$-stable parameters of the representations satisfy 
\begin{enumerate}
\item[{\rm{(a)}}]  $\varepsilon = \delta$, 
\item[{\rm{(b)}}]  $h_\pi \in \{ h_\Pi, h_\Pi -1\}$, 
\item[{\rm{(c)}}]  
$
\mu_0 \geq \nu_0 \geq \mu_1 \geq \dots \geq \mu_m \geq \nu_{m -1}\geq \mu_{m} = 0$.  
\end{enumerate}
\end{enumerate}
\end{theorem}

A detailed proof will be given in \cite{branching}.

\begin{remark}
A similar theorem was proved in  \cite[Thms. 4.1 and 4.2]{sbonvec}
 for the restriction  of  irreducible representations of $O(n+1,1)$ to $O(n,1)$ with
 trivial infinitesimal character $\rho$.  
\end{remark}

\

\noindent
{\bf Applications of the branching theorem}

\medskip
\noindent
{\bf \em Gross--Prasad conjectures.} \ 
The discussion in \cite[Chap.~13, Sect.~3.3]{sbonvec} shows
 that Theorem \ref{branch} in the special case
 where $h_{\Pi}=\left[\frac{n+1}{2}\right]$ and $h_{\pi}=\left[\frac{n}{2}\right]$ implies the following.  

\begin{theorem}
The { Gross--Prasad conjectures} are valid for all tempered representations with nonsingular infinitesimal character of $S O(n+1,1)$ and $S O(n,1)$,
 which are nontrivial on the center.  
\end{theorem}

\begin{remark}
(1)\enspace 
For tempered principal series representations $\Pi$ of $SO(n+1,1)$, 
 and $\pi$ of $SO(n,1)$
 which are nontrivial on each center, 
 it was proved in \cite[Chap.~11, Sect.~4]{sbonvec}.
For irreducible tempered representations 
 with trivial infinitesimal character $\rho$, 
 this was announced in \cite{sbonGP}
 and proved in \cite[Chap.~11, Sect.~5]{sbonvec}.

\noindent
(2)\enspace The branching law for nontempered representations 
 (Theorem \ref{branch})
 interpolates  between the classical branching laws of finite-dimensional representations \eqref{eqn:branchWeyl} and the branching laws of the conjecture by Gross and Prasad  for tempered representations.
\end{remark}

\medskip
\noindent
{\bf \em Periods.} \ 
For representations $\Pi$, $\pi$ of a real reductive Lie group $G$,
 respectively of a reductive subgroup $G'$,
 the space of symmetry breaking operators 
\[ 
     \operatorname{Hom}_{G'}(\Pi|_{G'}, \pi^{\vee}) 
\]
 and the space of $G' $-invariant continuous linear functionals
\[\operatorname{Hom}_{G'}(\Pi \otimes \pi,\bC)\]
 are naturally isomorphic to each other
 \cite[Thm. 5.4]{sbonvec}, 
 where $\pi^{\vee}$ denotes the contragredient representation of $\pi$
 in the category of admissible smooth representations
 (see Section \ref{subsec:III.2}), 
 and $\Pi \otimes \pi$ denotes the representation of $G'$
 acting on the outer tensor product representation
 $\Pi \boxtimes \pi$ of $G \times G'$ diagonally. 
Thus
 we may use symmetry breaking operators
 to determine $G'$-invariant continuous  linear functionals on $\Pi \otimes \pi$, {\it{i.e.}}, 
 periods. 
Hence, 
the branching theorem implies the following.  

\medskip 
 
\begin{theorem}
Suppose that the representations $\Pi $ and $\pi$ of $G$,
 respectively of $G'$ satisfy Assumption A
 $($see Section \ref{subsec:II.1.1}$)$
 with height $i$ respectively $j$. 
Then the following statements on the pair $(\Pi, \pi)$ are equivalent:
\begin{enumerate}
\item[{\rm{(i)}}] 
The representation $\Pi \boxtimes \pi$ has a nontrivial $G'$-period;
\item[{\rm{(ii)}}] $\Pi $ and $\pi$ have the same signature, 
$j=i$ or $i-1$ and their enhanced $\theta$-stable parameters satisfy the interlacing conditions of Theorem \ref{branch}.   
\end{enumerate}
\end{theorem}

\medskip
For representations with infinitesimal character $\rho$
 we proved in \cite[Chap.~10]{sbonvec} furthermore the following:

\begin{theorem}
\label{thm:testvector}
If $\Pi$ and $\pi$ are representations with trivial infinitesimal character $\rho$, 
 then any nonzero period
 does not vanish on the minimal $K$-type
 of the outer tensor product representation $\Pi \boxtimes \pi$.
\end{theorem}

\medskip 
In a special case we determine in \cite[Chap.~12]{sbonvec}
 the value of the period on vectors in the minimal $K$-type of a representation with infinitesimal character $\rho$. See also Theorem \ref{thm:period2}.  

\medskip

\noindent
{\bf \em Distinguished representations.}
Let $G$ be a reductive group and $H$ a reductive subgroup.  
We regard $H$ as a subgroup of the direct product group $G \times H$
 via the diagonal embedding $H \hookrightarrow G \times H$.  

\begin{definition}
Let $\psi$ be a one-dimensional representation of $H$. 
We say an admissible smooth representation  $\Pi$ of $G$ is
$(H,\psi)$-{\it{distinguished}}
 if 
\[ 
  (\operatorname{Hom}_{H}(\Pi \boxtimes \psi^{\vee}, {\mathbb{C}})
  \simeq)
  \operatorname{Hom}_{H}(\Pi|_H, \psi) \not = \{0\}.  
\]
If the character $\psi $ is trivial,
 we say $\Pi$ is $H$-{\it{distinguished}}.  
\end{definition}

Let $G=S O(n+1,1)$, 
 and ${\mathfrak{g}} \simeq {\mathfrak{s o}}(n+1,1)$
 its Lie algebra.  
We fix a fundamental Cartan subalgebra ${\mathfrak{h}}$ of ${\mathfrak {g}}$.  
For $0 \le i \le [\frac {n+1} 2]$, 
 there are $\theta$-stable parabolic subalgebras 
$
{\mathfrak {q}}_i\equiv {\mathfrak {q}}_i^+
$
$
=({\mathfrak {l}}_i)_{\mathbb{C}} +{\mathfrak {u}}_i
$
and 
$
{\mathfrak {q}}_i^-
=({\mathfrak {l}}_i)_{\mathbb{C}} +{\mathfrak {u}}_i^-
$
 in ${\mathfrak{g}}_{\mathbb{C}}=\operatorname{Lie}(G) \otimes_{\mathbb{R}} {\mathbb{C}}
={\mathfrak {u}}_i^- + ({\mathfrak {l}}_i)_{\mathbb{C}}+{\mathfrak {u}}_i
$
 such that 
 ${\mathfrak {q}}_i$ and ${\mathfrak {q}}_i^-$
 contain ${\mathfrak {h}}_{\mathbb{C}}$
 and that both the Levi subgroups $L_i$ of 
 ${\mathfrak {q}}={\mathfrak {q}}_i$ and ${\mathfrak{q}}_i^-$
 are given by 
$L_i \simeq S O(2)^i \times SO(n-2i+1,1)$, 
 see \cite[Lem.~14.38]{sbonvec}.

\medskip
There are two one-dimensional representations of $S O(n-2i+1,1)$
 when $n \ne 2i-1$:
 we denote by $\chi_+$ the trivial representation ${\bf{1}}$
 and by $\chi_-$ the nontrivial one.  
For $n=2i-1$, 
 we consider only $\chi_+$.  
We consider the differential $\lambda =(\lambda_1, \cdots, \lambda_i,0,\dots 0 ) $ of a one-dimensional representation of $L_i$ and assume that $\lambda$ satisfies the conditions in \cite[Chap.~14, Sect.~9]{sbonvec}
 and $\delta \in \{+,- \}$.  
We consider an irreducible one-dimensional $L_i$-module
\[
  {\mathbb{C}}_{\lambda} \boxtimes \chi_{\delta}
\]
and define an admissible smooth representation of $G$ of moderate growth
  denoted by $A_{\bq_i}(\lambda)_\delta$. 
 Its underlying $({\mathfrak{g}},K)$-module
 is given by the cohomological parabolic induction from 
 ${\mathfrak{q}}_i$, 
 see \cite{KV, VZ}. 
For $\delta=+$ we often omit the subscript.

\begin{theorem}
\label{thm:I.9}
Suppose that $\Pi \in {\mathcal A}$
 (see Section \ref{subsec:II.1.1} below for definition)
 is a representation of $S O(n+1,1)$ cohomologically induced 
 from a one-dimensional representation of a $\theta$-stable parabolic subalgebra ${\mathfrak{q}}_i$, i.e. that $\Pi= A_{\bq_i} (\lambda)$.
Then the height 
 (see Section \ref{sec:II.3})
 of $\Pi$ is $i$
 and $\Pi $ is $S O(n+1-i,1)$-distinguished.  
\end{theorem}

\begin{remark}
The proof of this theorem will be given 
 in a subsequent paper, 
 based on the work \cite{sbonvec}.  
For a different proof and perspective of this theorem,
 see \cite{K1}.

\end{remark} 

The irreducible representations
 with trivial infinitesimal character $\rho$ are obtained through cohomological induction. 
We set $\Pi_i :=A_{\frak{q}_i}(0)$.  
Then Theorem \ref{thm:I.9} may be regarded
 as a generalization of the following results in \cite[Thm.~12.4 and Lem.~15.10]{sbonvec}:
 
\begin{theorem}
\label{thm:period1}
Let $0 \le i \le n+1$.  
Then the representations $\Pi_{i}$ 
 of $G=SO(n+1,1)$ are  $SO(n+1-i,1)$-distinguished. 
\end{theorem}

For details see \cite[Chap.~12]{sbonvec}.

\medskip
\noindent
{\bf \em  A bilinear form on the $(\bg,K)$-cohomology.} \ 
In  \cite[Chap.~12, Sect.~3]{sbonvec}, 
 we considered the morphism on  $(\bg,K)$-cohomologies of representations  induced by a symmetry breaking operator:
Let $(G,G')=(SO(n+1,1),SO(n,1))$
 and $\Pi$, $\pi$ be irreducible  representations and $V, V'$ 
 irreducible finite-dimensional representations
 of $G$ and $G'$, 
respectively.  
By abuse of notation,
 we use the same symbols $\Pi$ and $\pi$
 to denote their underlying $({\mathfrak{g}},K)$-module
 and $({\mathfrak{g}}',K')$-module, 
 respectively,
 when we take $({\mathfrak{g}},K)$-cohomologies and $({\mathfrak{g}}',K')$-cohomologies.  
Suppose that 
there exists a symmetry breaking operator $T\colon \Pi \otimes V \to \pi \otimes V'$.  
Then the   symmetry breaking operator  $T \colon \Pi \otimes V \to \pi \otimes V'$  
induces for every $j$ a morphism
\[
   T^j \colon
   H^j({\mathfrak{g}}, K; \Pi \otimes V )
   \rightarrow
   H^{j}({\mathfrak{g}}', K'; \pi \otimes V')
\] 
on the $(\bg,K)$-cohomologies
 and a bilinear form 
\begin{align*}
   B_T \colon
   H^j({\mathfrak{g}}, K; \Pi \otimes V )
   \times 
   H^{n-j}({\mathfrak{g}}', K'; (\pi \otimes V')^{\vee}\otimes \chi_{(-1)^{n+1}})
   &\to {\mathbb{C}}
\end{align*}
where $(\pi \otimes V')^{\vee}$ denotes the contragredient representation
 of $\pi \otimes V'$. 
 
The induced morphism
 on the $(\bg,K)$-cohomologies may be zero,
 but in some special cases we conclude that it is nonzero.

\medskip
For $0 \le \ell \le [\frac{n+1}2]$
 and $\delta \in \{+, -\}$, 
 we denote by $\Pi_{\ell, \delta}$ the irreducible
 admissible smooth representation of $G=S O(n+1,1)$
 with underlying $({\mathfrak{g}},K)$-module
 $A_{\mathfrak{q}_{\ell}}(0)_{\delta}$
 (see Theorem \ref{thm:I.9} for notation).  
We also write simply $\Pi_{\ell}$
 for $\Pi_{\ell, \delta}$
 if $\delta=+$.  
We recall from \cite[Thm.~2.20]{sbonvec}
 that $\Pi_{\ell, \delta}$ is the unique submodule
 of the principal series representation
 $I_{\delta}(\Lambda^{\ell}({\mathbb{C}}^n, \ell))$
 for $\ell \ne \frac n 2$
  (see \eqref{eqn:IVlmd} below), 
 and is isomorphic to $I_{\delta}(V_+, \frac n 2) \simeq I_{\delta}(V_-, \frac n 2)$
 for $\ell = \frac n 2$
 where $V_+$ and $V_-$ are irreducible $SO(n)$-modules
 such that $\Lambda^{\frac n2}({\mathbb{C}}^n)=V_+ \oplus V_-$.  
We also recall from \cite[Prop.~15.11]{sbonvec}
 that the set of irreducible admissible representations
 of $G=S O(n+1,1)$ with the trivial ${\mathfrak{Z}}({\mathfrak{g}})$-infinitesimal character $\rho$ is classified as follows:
\begin{alignat*}{2}
&\{\Pi_{\ell, \delta}: 0 \le \ell \le \frac{n-1}2,\,\delta \in \{+,-\}\}
\cup \{\Pi_{\frac {n+1}2, +}\}
\quad
&&\text{if $n$ is odd, }
\\
&\{\Pi_{\ell, \delta}: 0 \le \ell \le \frac{n}2,\,\delta \in \{+,-\}\}
\quad
&&\text{if $n$ is even. }
\end{alignat*}
Analogous notation $\pi_{j,\varepsilon}$ is applied
 to the subgroup $G'=S O(n,1)$.

Then the following theorem follows from 
 \cite[Cor.~12.19 and Lem.~15.10]{sbonvec}
 for the nonvanishing of $B_T$
 and \cite[Thm.~15.19]{sbonvec} for the uniqueness of $T$.  

\begin{theorem}
Let $(G,G')=(SO(n+1,1),SO(n,1))$, 
 $0 \le i \le \frac n 2$, 
 and $\delta \in \{ \pm \}$.  
Let $T$ be a nontrivial symmetry breaking operator 
 $\Pi_{i,\delta} \to \pi_{i,\delta}$.  
\begin{enumerate}
\item[{\rm{(1)}}]
$T$ induces bilinear forms
\[
   B_T \colon 
   H^j({\mathfrak{g}}, K; \Pi_{i,\delta}) 
   \times 
   H^{n-j}({\mathfrak{g}}', K'; \pi_{i,(-1)^{n-1} \delta})
   \to {\mathbb{C}}
\quad
\text{for all $j$.}
\]
\item[{\rm{(2)}}]
The bilinear form $B_T$ is nonzero
 if and only if $j=i$ and $\delta=(-1)^i$.  
\end{enumerate}
\end{theorem}

In Section \ref{sec:6}, 
we also state a nonvanishing theorem for bilinear forms on the $(\bg,K)$-cohomologies of principal series representations, 
 see Theorem \ref{thm:Bnonzero}.

\medskip

\begin{remark}If $\Pi$ and $\pi$ have trivial infinitesimal character $\rho$ then $V$ and $V'$ are the trivial representations  
and the theorem
 follows from \cite[Chaps.~9 and 12]{sbonvec}.  
\end{remark}

Detailed proofs of the results 
 will be published elsewhere \cite{branching}.

\section{Representations with nonsingular integral infinitesimal character}

In this section
 we recall from \cite{sbonvec}
 the results  about principal series representations
 and irreducible representations of $G=S O(n+1,1)$
 and introduce their height and signature.

\subsection{Notation}

We are using the same notation and assumptions
 as in the book \cite{sbonvec}
 except that we do not use \lq\lq{bar}\rq\rq\
 for subgroups of $S O(n+1,1)$, 
 see Remark \ref{rem:I.1}. 
The proofs of most of the results stated  in this section can be found  in Chapter 2 and  Appendices I, II, and III therein.

We first recall some notation which is  up to small changes (see the comments in Remark \ref{rem:I.1})  the same as the notation 
 in the Memoir article \cite{sbon}.

Consider the standard quadratic form 
\begin{equation}
\label{eqn:quad}
      x_0^2 + x_1^2 +\dots +x_{n}^2-x_{n+1}^2
\end{equation}
 of signature $(n+1,1)$. 
We define $G$
 to be the indefinite special orthogonal group $SO(n+1,1)$
 that preserves the quadratic form 
 \eqref{eqn:quad} and the orientation.  
Let $G'$ be the stabilizer
 of the vector $e_{n}={}^t\! (0,0,\cdots,0,1,0)$.  
Then $G'$ is isomorphic to $SO(n,1)$.  
We set 
\begin{alignat}{3}
\label{eqn:K}
K &:=O(n+2) \cap G
&&=\{
\begin{pmatrix}
B &  
\\
  & \det B
\end{pmatrix}
:
B \in O(n+1) \} 
&& \simeq O(n+1), 
\\
K' &:= K \cap G' 
&&=\{
\begin{pmatrix}
B & & 
\\
  & 1 &
\\
  & & \det B
\end{pmatrix}
:
B \in O(n)
\}
&&\simeq O(n).  
\notag
\end{alignat}
Then $K$ and $K'$ are maximal compact subgroups
 of $G$ and $G'$, 
respectively.  

\subsection{Principal series representations}

Let ${\mathfrak {g}}={\mathfrak {so}}(n+1,1)$
 and ${\mathfrak {g}}'={\mathfrak {so}}(n,1)$
 be the Lie algebras of $G$
 and $G'$, 
 respectively.  
We take a hyperbolic element $H$
 as 
\begin{equation}
H :=
  E_{0,n+1} + E_{n+1,0} \in \mathfrak{g}', 
\end{equation}
and set
\begin{equation*}
{\mathfrak {a}}:={\mathbb{R}}H
\qquad
\text{ and }
A:=\exp {\mathfrak {a}}.  
\end{equation*}
Then the centralizers of $H$ in $G$ and $G'$
 are given by $M A$ and $M' A$, 
 respectively,
 where 
\begin{align*}
M :={}
  & \left\{
    \begin{pmatrix} 
        \varepsilon \\ & B \\ & & \varepsilon
    \end{pmatrix} :
    B \in SO(n), \   \varepsilon = \pm 1
    \right\} 
   & \simeq & SO(n) \times O(1), 
\\
M' :={}
  &  \left\{
    \begin{pmatrix} 
        \varepsilon \\ & B \\ & & 1 \\ & & & \varepsilon
    \end{pmatrix} :
    B \in SO(n-1):   \varepsilon = \pm 1
    \right\}
 &  \simeq & SO(n-1) \times O(1).  
\label{eqn:Mprime}
\end{align*}
We observe that 
 $\operatorname{ad}(H) \in \operatorname{End}_{\mathbb{R}}(\mathfrak{g})$
 has eigenvalues 
 $-1$, $0$, and $+1$.  
Let 
\[
{\mathfrak {g}} ={\mathfrak {n}}_- + ({\mathfrak {m}}+{\mathfrak {a}})+{\mathfrak {n}}_+
\]
be the corresponding eigenspace decomposition, 
 and $P$ a minimal parabolic subgroup 
 with Langlands decomposition $P=MAN_+$. 
Likewise, 
 $P':=M'AN_+'$ is a compatible Langlands decomposition
 of a minimal (also maximal) parabolic subgroup $P'$ of $G'$
 with Lie algebra 
\begin{equation}
\label{eqn:Lang}
{\mathfrak {p}}'={\mathfrak {m}}'+{\mathfrak {a}}+{\mathfrak {n}}_+'
 = ({\mathfrak {m}} \cap {\mathfrak {g}}')
  +({\mathfrak {a}} \cap {\mathfrak {g}}')
  +({\mathfrak {n}}_++{\mathfrak {g}}').  
\end{equation}

We note 
 that we have chosen $H \in {\mathfrak {g}}'$
 so that $P'=P \cap G'$
 and $A=\exp ({\mathbb{R}} H)$ is a common maximally split abelian subgroup
 in $P'$ and $P$.

\medskip
The character group of $O(1)$
 consists of two characters.  
We write $+$ for the trivial character ${\bf{1}}$,
 and $-$ for the nontrivial character.  
Since $M \simeq SO(n) \times O(1)$, 
 any irreducible representation
 of $M$ is the outer tensor product
 of an irreducible representation $(\sigma, V)$ of $S O(n)$
 and a character $\delta$ of $O(1)$.

Given $(\sigma,V) \in \widehat {SO(n)}$, 
 $\delta \in \{ \pm \} \simeq \widehat{O(1)}$, 
 and a character $e_{\lambda}(\exp (t H))=e^{\lambda t}$
 of $A$ for $\lambda \in {\mathbb{C}}$, 
 we define the (unnormalized) principal series representation
\begin{equation}
\label{eqn:IVlmd}
   I_{\delta}(V, \lambda)= \operatorname{Ind}_{P}^{G}
  (V \otimes \delta, \lambda) 
\end{equation}
of $G=SO(n+1,1)$ on the Fr{\'e}chet space 
 of smooth maps
 $f \colon G \to V$ 
 subject to 
\begin{multline*}
  f(g m m' e^{tH} n) 
  = 
  \sigma(m)^{-1} \delta(m') e^{-\lambda t} f(g)
\\
\text{for all
 $g \in G$, $m m' \in M \simeq SO(n) \times O(1)$, $t \in {\mathbb{R}}$, 
 $n \in N_+$.}
\end{multline*}

By a result of R.~Langlands \cite{La88}
 every irreducible nontempered representation with nonsingular integral infinitesimal character is isomorphic to the unique subrepresentation of a principal series representation $I_{\delta}(V, \lambda)$ with $\lambda < \frac n 2$.  
We denote it by  $\Pi_\delta(V,\lambda)$, 
 see also \cite[Chap.~15, Sect.~7]{sbonvec}.  

\medskip

\begin{definition}
We call the triple $(V,\delta ,\lambda)$ the {\it{Langlands parameter}}
 of the irreducible nontempered representation $\Pi_\delta (V,\lambda)$. 
\end{definition}

\subsubsection{The set $ {\mathcal A}$}
\label{subsec:II.1.1}

Since we highlight in the article representations which are of interest
 in number theory, 
 we consider from now on a subset of irreducible representations
 of special orthogonal groups.  
For results about irreducible representation of $SO(n+1,1)$
 in the general case see \cite[Chap.~15]{sbonvec}.

We start with irreducible finite-dimensional representations of $SO(n+1,1)$
 that are self-dual, 
 or equivalently,
 that are obtained as the restriction 
 of irreducible representations of $O(n+1,1)$. 
So for $n =2m$ we assume that the highest weight 
$(\mu_1,\dots, \mu_{m+1})$
 of an irreducible finite-dimensional representation $F$
 of $G=SO(2m+1,1)$
  is of the form $(\mu_1,\dots, \mu_m, 0)$.

\medskip
\noindent
{\bf Assumption A.}
Suppose that $\Pi$ is an irreducible representation
 of $G=SO(n+1,1)$ with regular integral infinitesimal character, 
 see \cite[Chap.~2, Sect.~1.4]{sbonvec}.  
We say that a representation $\Pi$ of $G$ satisfies {\it{Assumption}} A
 if it has the same infinitesimal character
 as a self-dual irreducible finite-dimensional representation of $G$.  
When $n$ is odd, 
 Assumption A is automatically satisfied.

\begin{notation} 
The set of irreducible representations
 of $G=SO(n+1,1)$
 satisfying Assumption A is denoted by $\mathcal A $.  
\end{notation}

For the convenience of the reader,
 we give a description of irreducible admissible representations
 of $G=SO(n+1,1)$
 satisfying Assumption A in Section \ref{subsec:Aeven}
 for $n$ even 
 and in Section \ref{subsec:Aodd} for $n$ odd.  

\subsubsection{Classification of the set ${\mathcal{A}}$ for $n$ even}
\label{subsec:Aeven}

Suppose $n=2m$.  
By using the highest weight,
 we write $V \in \widehat{SO(n)}$
 as $V=F^{SO(n)}(\sigma)$
 with $\sigma=(\sigma_1, \cdots, \sigma_m) \in {\mathbb{Z}}^m$
 satisfying $\sigma_1 \ge \cdots \ge \sigma_{m-1} \ge |\sigma_m|$.  
Then the contragredient representation of $V$
 is given as 
\[
   V^{\vee} \simeq F^{SO(n)}(\sigma_1, \cdots, \sigma_{m-1}, -\sigma_m).  
\]

Hence $V$ is self-dual if and only if $\sigma_m=0$.  

\begin{prop}
\label{prop:191313}
For $V=F^{SO(2m)}(\sigma)$, 
 we consider the following condition on $\lambda$:
\begin{equation}
\label{eqn:Redeven}
\lambda \in {\mathbb{Z}},\,\, 
\lambda < m,\,\, 
\lambda \not \in \{1-\sigma_1, \cdots, m-\sigma_m\}.  
\end{equation}
Then irreducible admissible representations of $SO(2m+1,1)$ in ${\mathcal{A}}$
 are classified as 
\begin{enumerate}
\item[$\bullet$]
{\rm{(nontempered case)}}
\[
\{\Pi_{\delta}(V,\lambda): \delta=\pm, \sigma_m=0, \text{ and $\lambda$ satisfies \eqref{eqn:Redeven}}\}, 
\]
or 
\item[$\bullet$]{\rm{(tempered case)}}
\[
\{I_{\delta}(V,m): \delta=\pm, \sigma_m>0\}.  
\]

\end{enumerate}
\end{prop}

We note that in the tempered case $V \not \simeq V^{\vee}$
 as $SO(n)$-modules
 because $\sigma_m >0$, 
 whereas there is a $G$-isomorphism, 
 see \cite[Prop.15.5]{sbonvec}:
\[
   I_{\delta}(V,m) \simeq I_{\delta}(V^{\vee},m).  
\]

\subsubsection{Classification of the set ${\mathcal{A}}$ for $n$ odd}
\label{subsec:Aodd}

Suppose $n=2m-1$.  
We write $V \in \widehat{SO(n)}$
 as $V=F^{SO(n)}(\sigma)$
 with $\sigma=(\sigma_1, \cdots, \sigma_{m-1}) \in {\mathbb{Z}}^{m-1}$
 satisfying $\sigma_1 \ge \cdots \ge \sigma_{m-1} \ge 0$.  

\begin{prop}
\label{prop:191315}
For $V=F^{SO(2m-1)}(\sigma)$, 
 we consider the following conditions on $\lambda$:
\begin{equation}
\label{eqn:Redodd}
\lambda \in {\mathbb{Z}},\,\, 
\lambda \le m-1,\,\, 
\text{ and }
\lambda \not \in \{1-\sigma_1, \cdots, m-1-\sigma_{m-1}\}.  
\end{equation}
\begin{equation}
\label{eqn:disclmd}
\lambda \in {\mathbb{Z}},\,\, 
m \le \lambda \le m-1+\sigma_{m-1}.  
\end{equation}
All irreducible admissible representations of $SO(2m,1)$ belong to ${\mathcal{A}}$, 
 which are classified as 
\begin{enumerate}
\item[$\bullet$]
{\rm{(nontempered case)}}
\[
\{\Pi_{\delta}(V,\lambda): \delta=\pm, \text{ $\lambda$ satisfies \eqref{eqn:Redodd}}\}, 
\]
or
\item[$\bullet$]{\rm{(discrete series)}}
\[
\{\Pi_{+}(V,\lambda): \sigma_{m-1}>0, \text{ $\lambda$ satisfies \eqref{eqn:disclmd}}\}.  
\]
\end{enumerate}
\end{prop}

We note that in the discrete series case 
 there is a $G$-isomorphism
\[
   \Pi_{+}(V,\lambda) \simeq \Pi_{-}(V,\lambda).  
\]

\subsection{The height of representations in $\mathcal A$}
\label{sec:II.3}

In this section,
 we define a {\it{height}}
\[
  h \colon {\mathcal{A}} \to \{0,1,\cdots,[\frac{n+1}2]\}
\]
for irreducible representations of $G=SO(n+1,1)$
 that belong to ${\mathcal{A}}$
 (see Section \ref{subsec:II.1.1}).

We recall that the $i$-th exterior representation
 $\Exterior^i({\mathbb{C}}^n)$
 of $S O(n)$ is irreducible
 if $2i \ne n$, 
 and splits into two irreducible representations
 if $2i =n$, 
 which may be written as 
\[
   \Exterior^{\frac n 2}({\mathbb{C}}^n)
   \simeq
   \Exterior^{\frac n 2}({\mathbb{C}}^n)^{(+)}
   \oplus
   \Exterior^{\frac n 2}({\mathbb{C}}^n)^{(-)}.  
\]
We define a finite family of principal series representations
 of $G=S O(n+1,1)$ by
\[
   I_{\delta}(i,i)
   :=
   \begin{cases}
   I_{\delta}(\Exterior^{i}({\mathbb{C}}^n),i)
   \quad
   &\text{for $i \ne \frac n 2$}, 
\\
   I_{\delta}(\Exterior^{i}({\mathbb{C}}^n)^{(\varepsilon)},i)
   \quad
   &\text{for $i = \frac n 2$}.  
   \end{cases}
\]
Then $I_{\delta}(i,i)$ has the trivial infinitesimal character $\rho$, 
 and it does not depend on the choice 
 $\varepsilon =+$ or $-$
 when $i=\frac n 2$, 
 see \cite[(15.5)]{sbonvec}.  

Suppose that  $\Pi$ is an irreducible nontempered representation of $SO(n+1,1)$
 in $\mathcal A$ with Langlands parameter $(V,\lambda,\delta)$.
The principal series representation $I_\delta(V, \lambda)$ can be obtained
 by using a translation functor from exactly one principal series representation 
$I_\delta(i,i) $, $i \in \{0, 1,\dots , [\frac{n+1}{2}] \}$
 without crossing a wall.  
See \cite[Thm.~16.24]{sbonvec}.

Following \cite[Chap.~14, Sect.~5 and Thm.~16.17]{sbonvec}, 
 we say that the principal series representation $I_\delta(V, \lambda)$
 has {\it{height}} $i$
 if it can be obtained from a principal series representations $I_\delta(i,i) $ by a translation functor without crossing walls.

\medskip
\begin{definition}
[height]
\begin{enumerate}
\item[{\rm{(1)}}]
Suppose that the representation $\Pi \in {\mathcal A}$ is not tempered and has Langlands parameter $(V,\lambda,\delta)$. 
If $I_\delta(V,\lambda)$ has \emph{height} $i$,
 we say that $\Pi$ has height $i$. 
\item[{\rm{(2)}}]
If $n=2m-1$
 we say that a discrete series representation $\Pi \in {\mathcal{A}}$
 has \emph{height} $m$.  
\item[{\rm{(3)}}]
If $n=2m$ and $\Pi \in {\mathcal A}$ is a tempered representation, 
 then we say that it has {\it{height}} $m$.
\end{enumerate}
\end{definition}

An explicit formula of the height can be derived from the case
 for $O(n+1,1)$ in \cite[Def.~14.26 and Chap.~15]{sbonvec}:
\begin{prop}
[height for nontempered representation]
\label{prop:height}
\begin{enumerate}
\item[{\rm{(1)}}]
Suppose $n=2m$.  
With notation as in Proposition \ref{prop:191313}, 
 the height $i$ of the nontempered representation 
 $\Pi_{\delta}(V,\lambda)$ takes the value
 in $\{0,1,\cdots,m-1\}$, 
 and is determined by the following inequalities:
\[
  i-\sigma_i < \lambda < i+1 -\sigma_{i+1}.  
\]
\item[{\rm{(2)}}]
Suppose $n=2m-1$.  
With notation as in Proposition \ref{prop:191315}, 
 the height of the nontempered representation 
 $\Pi_{\delta}(V,\lambda)$ takes the value
 in $\{0,1,\cdots,m-1\}$, 
 and is determined by the following condition:
\begin{enumerate}
\item[$\bullet$]
for $\lambda < m-1-\sigma_{m-1}$, 
 we  have $0 \le i \le m-2$ 
 and 
\[
  i-\sigma_i < \lambda < i+1 -\sigma_{i+1};  
\]
\item[$\bullet$]
for $m-1-\sigma_{m-1}<\lambda<m$, 
 we have $i=m-1$.   
\end{enumerate}
\end{enumerate}
\end{prop}

\begin{example}
For $\Pi \in {\mathcal{A}}$, 
 the height of $\Pi$ is zero 
 if and only if $\Pi$ is finite-dimensional.  
\end{example}

\subsection{Signatures of  representations in $\mathcal A$}
\label{subsec:II.4}
In this section,
 we define a {\it{signature}}
\[
\operatorname{sgn}\colon \mathcal{A} \to \{+, -, \pm \}
\]
for irreducible representations of $G=SO(n+1,1)$ 
 that belong to ${\mathcal{A}}$.  
We shall impose the condition 
\begin{equation}
\label{eqn:axiomsgn}
\operatorname{sgn}(\Pi \otimes \chi_-)=-\operatorname{sgn}\Pi.  
\end{equation}

\subsubsection{Tempered Representations}

We recall from Propositions \ref{prop:191313} and \ref{prop:191315}
 (see also \cite[Thms.~13.7 and 13.9]{sbonvec})
 that for every irreducible finite-dimensional representation $F$
 in $\mathcal A$
 there exist irreducible tempered representations
 with the same infinitesimal character.

Suppose $G=SO(2m+1,1)$ and $\lambda=m$ $(=\frac n2)$.  
The unitary principal series representation 
 $I_{\delta}(V,m)$ is tempered,
 and it has the same infinitesimal character
 as an irreducible finite-dimensional representation of $G$
 if and only if 
 $V$ is not self-dual.  
In this case,
 there is an isomorphism
 $I_{\delta}(V,m) \simeq I_{\delta}(V^{\vee},m)$
 as $SO(2m+1,1)$-modules
 as we saw in Section \ref{subsec:Aeven}.  
We define the signature of $I_{\delta}(V,m)$
 to be $\delta$.

\medskip
For $G=SO(2m,1)$, 
 there is exactly one discrete series representation $\Pi$ of $G$
 having the same infinitesimal character as $F$
 (see Proposition \ref{prop:191315}). 
Furthermore there is an isomorphism
\[ \Pi \otimes \chi_- \simeq \Pi \]
 as $G$-modules, 
 hence we define the signature
 of a discrete series representation $\Pi$ to be $\pm$.

\subsubsection{Nontempered representations}
\label{subsec:II.2.2}
An irreducible representation $\Pi$ in $\mathcal A$
 of $G=SO(n+1,1)$,
 which is not tempered, is isomorphic to a subrepresentation
 of a principal series representation  $I_{\delta}(V, \lambda)$
 of height $i \in \{0,1,\cdots,[\frac {n-1}{2}]\}$, 
 see Proposition \ref{prop:height}.  
It has a Langlands parameter  $(V,\delta,\lambda)$ where $V$ is a self-dual representation of $SO(n)$, $\lambda < \frac n2$ and $\delta \in \{ +,- \}$.    
We refer to $\delta (-1)^{i - \lambda}$ as the {signature of the representation $\Pi$.  }
Since the Langlands parameter of a nontempered representation is unique, the signature is unique as far as $\Pi$ is not tempered.

\begin{remark}
If a  representation $\Pi$ has signature $\delta$, 
 then one sees $\Pi \otimes \chi_-$ has signature $-\delta$.
\end{remark}

\begin{example} 
[one-dimensional representations]
\label{ex:chi}
There are two one-dimensional representations of 
$G=SO(n+1,1)$: we denote the trivial representation ${\bf{1}}$ by $\chi_+$
and the nontrivial one by $\chi_-$. 
The Langlands parameter of $\chi_+$ and $\chi_-$ is 
 $(+, {\bf{1}}, 0)$ and $(-, {\bf{1}}, 0)$, 
 respectively,
 and their height is 0.  
Hence the representation $\chi_{+}$ has signature $+$
 and the representation $\chi_{-}$ has signature $-$. 
\end{example}

\begin{example}
[irreducible representations
 with trivial infinitesimal character $\rho$]
\label{ex:reprho}
If $V$ is the representation of $SO(n)$
 on the $i$-th exterior tensor space $\bigwedge^i({\mathbb{C}}^{n})$
 ($2i \ne n$), 
 we write for simplicity  $I_{\delta}(i, \lambda)$
 instead of  $I_{\delta}(V, \lambda)$.  
Then the $SO(n)$-isomorphism
 on the exterior representations
$
  \bigwedge^i({\mathbb{C}}^{n}) \simeq \bigwedge^{n-i}({\mathbb{C}}^{n})
$
 leads us to the following $G$-isomorphism:
\[
   I_{\delta}(i, \lambda) \simeq I_{\delta}(n-i, \lambda).  
\]
If $n$ is even and $n=2i$, 
the exterior representation $\bigwedge^i({\mathbb{C}}^{n})$
 splits into two irreducible representations
 of $SO(n)$:
\[
{\Exterior}^{\frac n 2}({\mathbb{C}}^{n})
\simeq
{\bigwedge}^{\frac n 2}({\mathbb{C}}^{n})_+
\oplus
{\bigwedge}^{\frac n 2}({\mathbb{C}}^{n})_-
\]
with highest weights $(1,\cdots,1,1)$ and $(1,\cdots,1,-1)$, 
 respectively, 
 with respect to a fixed positive system for ${\mathfrak {so}}(n,{\mathbb{C}})$.  
Accordingly,
 we have a direct sum decomposition
 of the induced representation:
\begin{align}
   \operatorname{Ind}_{P}^{G}
  ({\Exterior}^{\frac n 2}({\mathbb{C}}^{n}) \otimes \delta, \lambda)
&=
  I_{\delta}({\bigwedge}^{\frac n 2}({\mathbb{C}}^{n})_+, \lambda)
\oplus
   I_{\delta}({\bigwedge}^{\frac n 2}({\mathbb{C}}^{n})_-, \lambda), 
\notag
\intertext{which we shall write as}
  I_{\delta}\left( \frac n 2, \lambda \right)
  &=
  I_{\delta}^{(+)} \left( \frac n 2, \lambda \right)
\oplus
   I_{\delta}^{(-)} \left( \frac n 2, \lambda \right) .  
\label{eqn:split}
\end{align}
The representations
$I_{\delta}^{(+)} \left( \frac n 2, \lambda \right)$ and  $I_{\delta}^{(-)} \left( \frac n 2, \lambda \right)$ are isomorphic to each other
 and have the signature
$ \delta$.
\end{example}

Let ${\mathfrak {Z}}({\mathfrak {g}})$ be the center of the enveloping algebra
 $U({\mathfrak {g}})$ 
 of the complexified Lie algebra 
 ${\mathfrak {g}}_{\mathbb{C}}={\mathfrak {g}} \otimes_{\mathbb{R}}{\mathbb{C}}
 \simeq {\mathfrak {s o}}(n+2, {\mathbb{C}})$.  
Via the Harish-Chandra isomorphism, 
 the ${\mathfrak {Z}}({\mathfrak {g}})$-infinitesimal character
 of the trivial one-dimensional representation ${\bf{1}}$
 is given by 
\[
\rho= \left( \frac{n}{2}, \frac{n}{2}-1,\cdots,\frac{n}{2}-[\frac{n}{2}]
      \right)
\]
 in the standard coordinates
 of the Cartan subalgebra of
$
   {\mathfrak {g}}_{\mathbb{C}}
$
$
   =
$
$
   {\mathfrak {so}}(n+2,{\mathbb{C}})
$, 
 whereas up to conjugation by the Weyl group the infinitesimal character of 
 $I_{\delta}^{(\pm)}(i, \lambda)$
 (when $2i \le n$)
 is given by 
\begin{equation}
\label{eqn:Zginf}
  \left(
   \frac n 2, \frac n 2-1, \cdots, \frac n 2-i+1, 
   \widehat{\frac n 2 -i}, 
   \frac n 2-i-1, 
   \cdots, 
   \frac n 2 - [\frac n 2], 
   \lambda-\frac n 2
  \right).  
\end{equation}

The irreducible representation $\Pi_{i,\delta}$ has height $\operatorname{min}(i,n-i)$.  
Since $\lambda =i$, 
 the signature is equal to $\delta$
 if $2i \le n$
 and to $\delta(-1)^n$ if $2i \ge n$.
The irreducible tempered representations are denoted by $\Pi_{m}$
 if $n=2m-1$ and $\Pi_{m,\delta}$ if $n= 2m$.
See \cite[Chap.~2, Sect.~4.5]{sbonvec} for $O(n+1,1)$
 and \cite[Chap.~15, Sect.~5]{sbonvec} for $S O(n+1,1)$ in detail.

For the group $G'=SO(n,1)$, 
 we shall use the notation $J_{\varepsilon}(j,\nu)$
 for the unnormalized parabolic induction
$
   \operatorname{Ind}_{P'}^{G'}
   (\bigwedge^j({\mathbb{C}}^{n-1})\otimes \varepsilon,\nu)
$
 for $0 \le j \le n-1$, 
 $\varepsilon \in \{+,-\}$, 
 and $\nu \in {\mathbb{C}}$.  
The irreducible representations  are denoted by $\pi _{j,\varepsilon}$ respectively $\pi_j$.

\subsubsection{Hasse and Standard sequences}

The notion of the height of representation in $\mathcal A$
 is motivated by the Hasse sequences in \cite[Chap.~13]{sbonvec}, 
 which were defined for the full orthogonal group 
 $O(n+1,1)$. 
We adapt the definition 
 for the special orthogonal group $G=S O(n+1,1)$
 as follows.  
Let $n=2m$ or $2m-1$.  
For every irreducible finite-dimensional representation $F$
 of the group $G$,
 there exists a unique sequence 
\begin{center}
\begin{tabular}{ccccccccccc}
& &$U_0$&, &\dots   & , & $U_{m-1} $ &, & $U_{m}$ 
\end{tabular}
\end{center}
of irreducible admissible smooth representations 
 $U_i\equiv U_i(F)$ of $G$ such that 
\begin{enumerate}
\item  $U_0 \simeq F$; 
\item  consecutive representations are composition factors of a principal series representation;
\item
$U_i$ $(0 \le i \le m)$ are pairwise inequivalent as $G$-modules.  
\end{enumerate}

\medskip

\begin{definition}
[Hasse sequence and standard sequence]
We refer to the sequence 
\begin{center}
\begin{tabular}{ccccccccccc}
& &$U_0$&, &\dots   & , & $U_{m-1} $ &, & $U_{m}$ 
\end{tabular}
\end{center}
 as the 
{\it {Hasse sequence}}
 of irreducible representations
 starting with the finite-dimensional representation $U_0=F$.  
We shall write $U_j(F)$ for $U_j$
 if we emphasize the sequence $\{U_j(F)\}$
 starts with $U_0=F$, 
 and  we refer to 
\begin{center}
\begin{tabular}{ccccccccccc}
& ${\Pi}_0:=U_0$&, &\dots   & , 
& ${\Pi}_{m-1}:=U_{m-1}\otimes (\chi_{-})^{m-1} $ &, 
& $\Pi_m:=U_{m}\otimes (\chi_{-})^m$ 
\end{tabular}
\end{center}
as the
 {\it{standard sequence}}
 of irreducible representations ${\Pi}_i={\Pi}_i(F)$
 starting with ${\Pi}_0=U_0=F$, 
 where $\chi_{-}$  
 is the nontrivial  one-dimensional representation of $G$
 defined in Example \ref{ex:chi}. 
\end{definition}

\medskip
More details about the standard sequence for $O(n+1,1)$ can be found in \cite[Chap.~13]{sbonvec}, from 
which the case for the normal subgroup $G=S O(n+1,1)$ is derived as follows.  
\begin{theorem}
Let $G=S O(n+1,1)$.  
\par\noindent
{\rm{(1)}}\enspace
Suppose that 
 $\Pi \in {\mathcal A}$ is not a discrete series representation
 of $G$,
 having height $j$ and signature $\delta$. 
Then there exists exactly one 
 irreducible finite-dimensional representation $F$ of $G$
 with signature $\delta$ so that 
$\Pi $ is the $j$-th representation in the standard sequence starting with $F$.  
\par\noindent
{\rm{(2)}}\enspace
Suppose that  $\Pi \in {\mathcal A}$ is a discrete series representation with signature $\delta \in \{\pm\}$.  
Then $n$ is odd, 
 and there exists a unique irreducible finite-dimensional representation $F$
 of $G$
 with signature $+$,
 so that  $\Pi$ is the  $\frac{n+1}{2}$-th representation in the standard sequence starting with $F$ and with $F \otimes \chi_-$.
\end{theorem}

\begin{remark}  
The  last representation in the standard sequence starting
 at an irreducible finite-dimensional representation $F \in {\mathcal A}$ is tempered.
\end{remark}

\subsubsection{The $\theta$-stable and enhanced $\theta $-stable parameters of irreducible representations}

We summarize the results in \cite[Chap.~14]{sbonvec}
 and give a parametrization of irreducible subquotients
 of the principal series representations
  $ I_{\delta}(V,\lambda)$
 of the group $G=SO(n+1,1)$ in terms of cohomological parabolic induction.  

We recall quickly 
 cohomological parabolic induction.  
A basic reference is Vogan \cite{Vogan81} and Knapp--Vogan \cite{KV}.  
We begin with a {\it{connected}} real reductive Lie group $G$.  
Let $K$ be a maximal compact subgroup, 
 and $\theta$ the corresponding Cartan involution.  
Given an element $X \in {\mathfrak{k}}$, 
 the complexified Lie algebra
\[
   {\mathfrak{g}}_{\mathbb{C}}={\operatorname{Lie}}(G) \otimes_{\mathbb{R}}
{\mathbb{C}}
\]
 is decomposed into the eigenspaces
 of $\sqrt{-1}\operatorname{ad}(X)$, 
 and we write 
\[
   {\mathfrak{g}}_{\mathbb{C}}
   ={\mathfrak{u}}_- + {\mathfrak{l}}_{\mathbb{C}} + {\mathfrak{u}}
\]
 for the sum of the eigenspaces 
 with negative, zero, 
 and positive eigenvalues.  
Then ${\mathfrak{q}}:={\mathfrak{l}}_{\mathbb{C}}+{\mathfrak{u}}$
 is a $\theta$-stable parabolic subalgebra 
 with Levi subgroup 
\begin{equation}
\label{eqn:LeviLq}
   L =\{g \in G: {\operatorname{Ad}}(g) {\mathfrak{q}}={\mathfrak{q}}\}.  
\end{equation}
The homogeneous space $G/L$ is endowed
 with a $G$-invariant complex manifold structure 
 such that its holomorphic cotangent bundle is given
 as $G \times_L {\mathfrak{u}}$.  
As an algebraic analogue of Dolbeault cohomology groups
 for $G$-equivariant holomorphic vector bundle over $G/L$, 
 Zuckerman introduced a cohomological parabolic induction functor
 ${\mathcal{R}}_{\mathfrak{q}}^j(\cdot \otimes {\mathbb{C}}
_{\rho({\mathfrak{u}})})$ ($j \in {\mathbb{N}}$) from the category
 of $({\mathfrak{l}}, L \cap K)$-modules
 to the category of $({\mathfrak{g}}, K)$-modules. 
We adopt here the normalization of the cohomological parabolic induction
 ${\mathcal{R}}_{{\mathfrak {q}}}^{j}$ from a $\theta$-stable parabolic subalgebra
 ${\mathfrak {q}}={\mathfrak {l}}_{\mathbb{C}}+{\mathfrak {u}}$
 so that the ${\mathfrak{Z}}({\mathfrak {g}})$-infinitesimal character
 of the $({\mathfrak{g}},K)$-module 
 ${\mathcal{R}}_{{\mathfrak {q}}}^{j}(F)$ equals
\[
\text{
 the ${\mathfrak{Z}}({\mathfrak {l}})$-infinitesimal character
 of the ${\mathfrak {l}}$-module $F$
}
\]
 modulo the Weyl group 
 via the Harish-Chandra isomorphism.  

\medskip

For each $i$ with $0 \le i \le [\frac {n+1} 2]$, 
 there are $\theta$-stable parabolic subalgebras 
$
{\mathfrak {q}}_i\equiv {\mathfrak {q}}_i^+
$
$
=({\mathfrak {l}}_i)_{\mathbb{C}} +{\mathfrak {u}}_i
$
and 
$
{\mathfrak {q}}_i^-
=({\mathfrak {l}}_i)_{\mathbb{C}} +{\mathfrak {u}}_i^-
$
 in ${\mathfrak {g}}_{\mathbb{C}}=\operatorname{Lie}(G) \otimes_{\mathbb{R}} {\mathbb{C}}$
such that ${\mathfrak {q}}_i$ and ${\mathfrak {q}}_i^-$
 contain a fundamental Cartan subalgebra 
 ${\mathfrak {h}}$.  
The Levi subgroup $L_i = N_G({\mathfrak {q}}_i)$ 
 of the $\theta$-stable parabolic subalgebra 
 ${\mathfrak {q}}_i$ and ${\mathfrak{q}}_i^-$
 is isomorphic to  
$L_i=SO(2)^i \times SO(n-2i+1,1)$.

We set 
\[
  \Lambda^+(N):=
  \{(\lambda_1, \cdots, \lambda_N)\in {\mathbb{Z}}^N:
    \lambda_1 \ge \lambda_2 \ge \cdots \ge \lambda_N \ge 0\}.  
\]
For $\nu=(\nu_1, \cdots, \nu_i) \in {\mathbb{Z}}^i$, 
 $\mu \in \Lambda^+([\frac n 2]-i+1)$, 
 and $\delta \in \{+,- \}$, 
 we consider an irreducible finite-dimensional $L_i$-module
\[
  F^{O(n-2i+1,1)}(\mu)_{\delta} \otimes {\mathbb{C}}_{\nu}
\]
and define an admissible smooth representation of $G$ of moderate growth, 
 whose underlying $({\mathfrak{g}},K)$-module
 is given by the cohomological parabolic induction
\begin{equation}
  {\mathcal{R}}_{{\mathfrak {q}}_i}^{S_i}(F^{S O(n-2i+1,1)}(\mu)_{\delta} \otimes {\mathbb{C}}_{\nu+\rho({\mathfrak{u}}_i)})
\end{equation}
 of degree $S_i$,  
where we set
\begin{equation}
   S_i:= \dim_{\mathbb{C}} ({\mathfrak {u}}_i \cap {\mathfrak {k}}_{\mathbb{C}})       =i(n-i).  
\end{equation}
It is  denoted by 
\[
 ( {\nu_1, \cdots, \nu_i \ || \  \mu_1, \cdots, \mu_{[\frac n 2]-i+1})_\delta}.
\]  
We note
 that if $i=0$
 then $( || \ {\mu_1, \cdots, \mu_{[\frac n 2]+1}})_{\delta}$ is
 finite-dimensional.  

\begin{definition}
We call $({\nu_1, \cdots, \nu_i \ ||\ \mu_1, \cdots, \mu_{[\frac n 2]-i+1})_
{\delta}}$
 the 
{\bf $\theta$-stable parameter }
 of the representation $ {\mathcal{R}}_{{\mathfrak {q}}_i}^{S_i}(F^{O(n-2i+1,1)}(\mu)_{\delta} \otimes {\mathbb{C}}_{\nu+\rho({\mathfrak{u}}_i)})$.  
\end{definition}

\begin{remark}
By \cite[Chaps.~14 and 16]{sbonvec},
 the double bars $||$ 
 in the $\theta$-stable parameter of a representation
 in $\mathcal A$
 of height $i$ 
 are before the $i+1$-th entry.
\end{remark}

\begin{remark}
In the introduction we refer to $(\mu, i,\delta)$ as the {\bf enhanced $\theta$-stable parameter }of the representation with $\theta$-stable parameter 
$({\mu_1, \cdots, \mu_i \ ||\ \mu_{i+1}, \cdots, \mu_{[\frac n 2]+1})_
{\delta}}$
\end{remark}

\begin{example}
\label{ex:II.14}
\begin{enumerate} 
\item[{\rm{(1)}}] 
An irreducible finite-dimensional  representation $F^G(\mu)_\delta$ has the $\theta$-stable parameter
$(|| \mu_1,\mu_2, \dots, \mu_{[\frac n 2]+1})_\delta $.
\item[{\rm{(2)}}] 
The $\theta$-stable parameter of a representation
 of height $i$ with trivial infinitesimal character $\rho$ is 
\[(0,0,\dots, 0\ || \ 0, \dots, 0)_\delta \]
where the  double bars  $||$ are before the $i+1$-th zero 
(see \cite[Chap.~14, Sect.~9.3]{sbonvec}).  
\item[{\rm{(3)}}] 
The  representations $\Pi$ in $\mathcal A$ with $\theta$-stable parameter $(\lambda_1,\lambda_2, \dots, \lambda_i\ || \ 0, \dots, 0)_\delta$ are unitary and are often referred to as $A_\bq (\lambda)_\delta$.
There exists a finite-dimensional representation $V$ of $G$
 so that $H^*(\bg,K;\Pi\otimes V) \not = \{0\}$, 
 see \cite{VZ}.  
\end{enumerate}
\end{example}

\subsubsection{The Hasse and Standard sequences in $\theta$-stable parameters.}

We set $m:=[\frac{n+1}{2}]$, 
 namely $n=2m-1$ or $2m$.  
Let $F=F^G(s_0,\dots , s_{[\frac n 2]})_\delta$ be an irreducible finite-dimensional representation of $G=SO(n+1,1)$,
 and $U_i \equiv U_i(F)$
 ($0 \le i \le [\frac{n+1}{2}]$) be the Hasse sequence with $U_0 \simeq F$.  
In \cite[Chap.~14]{sbonvec} we show:

\begin{theorem}
Let $n=2m$
 and $0 \le i \le m$.  
\begin{enumerate}
\item[{\rm{(1)}}]
{\rm{(Hasse sequence)}}\enspace
$U_i(F) \simeq  (s_0, \cdots, s_{i-1}||s_i, \cdots, s_m)_
{ (-1)^{i-s_i} \delta}$.  
\item[{\rm{(2)}}]
{\rm{(standard sequence)}}\enspace
$U_i(F) \otimes \chi_{-}^i \simeq 
 (s_0, \cdots, s_{i-1} || \ s_i, \cdots, s_m)_{(-1)^{s_i}\delta}$. 
\end{enumerate} 
\end{theorem}

\medskip

\noindent
The case $n$ odd is given similarly as follows.  
\begin{theorem}
Let $n=2m-1$, 
 and $0 \le i \le m-1$.  
\begin{enumerate}
\item[{\rm{(1)}}]
{\rm{(Hasse sequence)}}\enspace
$U_i(F) \simeq  (s_0, \cdots, s_{i-1} || \ s_i, \cdots, s_{m-1})_
{(-1)^{i-s_i} \delta}$.  
\item[{\rm{(2)}}]
{\rm{(standard sequence)}}\enspace
$U_i(F) \otimes \chi_{-}^i \simeq 
\ (s_0, \cdots, s_{i-1}|| \ s_i, \cdots, s_{m-1})_{ (-1)^{s_i}\delta}$. 
\end{enumerate}
\end{theorem}

\medskip

\section{The restriction of  representations of $SO(n+1,1)$ in $\mathcal A$  to the subgroup $SO(n,1)$}
In this section we discuss the branching law for  the restriction of irreducible representations $\Pi \in {\mathcal A}$
 of $SO(n+1,1)$ to the subgroup $SO(n,1)$. 
We state it for infinite-dimensional representations in Langlands parameter
 and $\theta$-stable parameters as well in the language of height and signature.A branching law for irreducible representations without the assumption $\Pi \in {\mathcal A}$ will appear in \cite{branching}.

\medskip

\subsection{Branching laws for finite-dimensional representations}
We first recall the branching laws for finite-dimensional representations.  As in the classical branching law for $SO(N) \downarrow SO(N-1)$ 
  the irreducible decomposition of finite-dimensional representations
 of $SO(N,1)$
 when restricted to the subgroup $SO(N-1,1)$
is  as follows:
\begin{theorem}
[branching rule for $S O(N,1) \downarrow S O(N-1,1)$]
Let $N \geq 2$.  
Suppose that $(\lambda_1, \cdots,\lambda_{[\frac{N+1}2]})
 \in \Lambda^+([\frac{N+1}2])$
 and $\delta \in \{+, -\}$.  
Then the irreducible finite-dimensional representation
 $ F^{O(N,1)}(\lambda_1, \cdots,\lambda_{[\frac{N+1}2]})_{\delta}$
  of $SO(N,1)$ decomposes
 into a multiplicity-free sum
 of irreducible representations
 of $SO(N-1,1)$ as follows:
\begin{equation*}
   F^{SO(N,1)}(\lambda_1, \cdots,\lambda_{[\frac{N+1}2]})_{\delta}|_{SO(N-1,1)}
   \simeq
  \bigoplus
   F^{SO(N-1,1)}(\nu_1, \cdots, \nu_{[\frac{N}{2}]})_{\delta}, 
\end{equation*}
where the summation is taken over
 $(\nu_1, \cdots, \nu_{[\frac{N}{2}]}) \in {\mathbb{Z}}^{[\frac N 2]}$
 subject to 
\begin{alignat*}{2}
&\lambda_{1} \ge \nu_{1} \ge \lambda_{2} \ge \cdots \ge \nu_{\frac N 2} \ge 0
\quad
&&\text{for $N$ even}, 
\\
&\lambda_{1} \ge \nu_{1} \ge \lambda_{2} \ge \cdots \ge \nu_{\frac {N-1} 2}
\ge \lambda_{\frac {N+1} 2}
\quad
&&\text{for $N$ odd}.   
\end{alignat*} 
\end{theorem}

\medskip

\begin{example}
We proved the branching rule for the restriction of the representations of $SO(n,1)$ on the space of harmonic polynomials  in  \cite[Prop.~2.3]{sbon}.
\end{example}

\subsection{Symmetry breaking operators}
\label{subsec:III.2}
Irreducible infinite-dimensional representations of $G$
 typically do not decompose
 into a direct sum of irreducible representations of $G$
 when restricted to a noncompact subgroup $G'$, 
 see \cite{xkInvent98} for details. 
To obtain information about the restriction and the branching laws we have to proceed differently.

For a continuous representation $\Pi$ of $G$
 on a complete, 
 locally convex topological vector space ${\mathcal{H}}$,
 the space ${\mathcal{H}}^\infty $ of $C^\infty$-vectors of ${\mathcal{H}}$
 is naturally endowed with a Fr{\'e}chet topology,
 and $(\Pi,{\mathcal{H}})$ induces a continuous representation $\Pi^{\infty}$ of $G$
 on ${\mathcal{H}}^\infty$.  
If $\Pi$ is an admissible representation
 of finite length on a Banach space ${\mathcal{H}}$, 
 then the Fr{\'e}chet representation 
 $(\Pi^{\infty}, {\mathcal{H}}^{\infty})$, 
 which we refer to as an 
 {\it{admissible smooth representation}}, 
 depends only 
 on the underlying $({\mathfrak {g}}, K)$-module
 ${\mathcal{H}}_K$.  
In the context of asymptotic behaviour 
 of matrix coefficients, 
 these representations
 are also referred to as an admissible representations
 of moderate growth 
 \cite[Chap.~11]{W}.  
We shall work with these representations
 and write simply $\Pi$ for $\Pi^{\infty}$.  
We denote by ${\operatorname{Irr}}(G)$
 the set of equivalence classes
 of irreducible admissible smooth representations.  
We also sometimes call these representations
 \lq\lq{irreducible admissible representations}\rq\rq\
 for simplicity.  

Given another admissible smooth representation $\pi$
 of a reductive subgroup $G'$, 
 we consider the space of continuous $G'$-intertwining operators 
 ({\it{symmetry breaking operators}})
\[ \operatorname{Hom}_{G'} ({\Pi}|_{G'}, {\pi}) .\] 
If $G=G'$ then these operators include the Knapp--Stein operators \cite{KS}
 and the differential intertwining operators 
 studied by B.~Kostant \cite{Kos}. 
Including the general case where $G \not = G'$, 
 we define now the multiplicity of $\pi$
 occurring in the restriction $\Pi|_{G'}$ as follows.  

\begin{definition}
[multiplicity]
For $G \supset G'$, 
 we say
\[
   m(\Pi,\pi) :=\operatorname{dim}_{\mathbb{C}} \operatorname{Hom}_{G'} ({\Pi}|_{G'}, {\pi}) 
\]
the  {\it{multiplicity}} of $\pi$
 occurring in the restriction $\Pi|_{G'}$.
\end{definition}

A finiteness criterion and a uniformly boundedness criterion
 are proved in \cite[Thms.~C and D]{xko}.  
Moreover, 
 by a result of B.~Sun and C.~B.~Zhu \cite{SunZhu}, 
 our assumptions imply that the multiplicities are either 0 or 1. 
These multiplicities yield important information of the restriction of $\Pi $ to $G'$,
 as we will see in the applications in the next part of the article.

\medskip
\subsection{Branching laws for representations in $\mathcal A$ : \\ First   formulation}

\begin{theorem}
\label{thm:branch1}
Let $F$ be an irreducible finite-dimensional representations
 of $G=SO(n+1,1)$,
 and $\{\Pi_i(F)\}$ be the standard sequence
 starting at $\Pi_0(F)=F$.  
Let $F'$ be an irreducible finite-dimensional representation of the subgroup $G'=SO(n,1)$, 
 and $\{\pi_j(F')\}$ the standard sequence starting 
 at $\pi_0(F')=F'$.  
Assume that 
\[ 
   {\operatorname{Hom}}_{G'}(F|_{G'},F') \not = \{0\}.  
\]
Then  symmetry breaking  for the representations $\Pi_i(F)$, $\pi_j(F')$
 in the standard sequences  is presented graphically
 in Diagrams \ref{tab:SBodd} and \ref{tab:SBeven}. 
In the first row are representations of G, in the second row are representations of $G'$. 
Nontrivial symmetry breaking operators are represented by arrows, 
 namely,
 there exist nonzero symmetry breaking operators between 2 representations 
 if and only if there are arrows
 in the Diagrams \ref{tab:SBodd} and \ref{tab:SBeven}.    
\end{theorem}

\begin{figure}[htp]
\color{black}
\caption{Symmetry breaking for $SO(2m+1,1) \downarrow SO(2m,1)$}
\begin{center}
\begin{tabular}{@{}c@{~}c@{~}c@{~}c@{~}c@{~}c@{~}c@{~}c@{~}c@{~}c@{}}
$\Pi_0(F)$
& 
&$\Pi_1(F)$
&  
&\dots
&
&$\Pi_{m-1}(F)$
& 
&$\Pi_{m}(F)$ 
\\
$\downarrow$ 
&$\swarrow$
& $\downarrow$
& $\swarrow$ 
& 
& $\swarrow$ 
& $ \downarrow $
&  $\swarrow $  
&  $\downarrow$ 
\\
$\pi_0(F')$& &$\pi_1(F')$
& 
&\dots 
&
& $\pi_{m-1}(F')$ 
& 
& $\pi_{m}(F')$ 
\end{tabular}
\end{center}
\label{tab:SBodd}
\end{figure}%

\medskip
\begin{figure}[htp]
\color{black}
\caption{Symmetry breaking for $SO(2m+2,1) \downarrow SO(2m+1,1)$ }
\begin{center}
\begin{tabular}{@{}c@{~}c@{~}c@{~}c@{~}c@{~}c@{~}c@{~}c@{~}c@{~}c@{~}c@{}}
$\Pi_0(F)$& &$\Pi_1(F)$ & & \dots &  & $\Pi_{m-1}(F)$& & $\Pi_{m}(F)$ & & $\Pi_{m+1}(F)$\\
$\downarrow$ &$\swarrow$& $\downarrow $& $\swarrow$ &  & $\swarrow$ & $ \downarrow $& $\swarrow $ & $\downarrow$ & $\swarrow$\\
$\pi_0(F')$& &$\pi_1(F')$& &\dots  & & $\pi_{m-1}(F')$ & & $\pi_{m}(F')$ 
\end{tabular}
\end{center}
\label{tab:SBeven}
\end{figure}%

\subsection{ Branching laws for representations in $\mathcal A$:\\ Second formulation}
Let $F^ G( \mu)_{\delta}$ and $F^{G'}( \nu )_{\delta}$ be  irreducible finite-dimensional representations in $\mathcal A$ of $G=SO(n+1,1)$, 
 respectively  of the subgroup $G'=SO(n,1)$, 
 where $\mu \in \Lambda^+([\frac{n+2}{2}])$, 
 $\nu \in \Lambda^+([\frac{n+1}{2}])$, 
 and $\delta\in \{+,-\}$.

Suppose that 
\[
  {\operatorname{Hom}}_{G'} (F^{G}{(\mu})_{\delta}|_{G'}, F^{G'}{(\nu )}_{\delta})\not = \{0\}.
\]

If $n=2m$, 
then $\mu =(\mu_0, \cdots, \mu_{m+1})\in \Lambda^+(m+1)$
 and $\nu =(\nu_0, \cdots, \nu_{m})\in \Lambda^+(m)$ and 
 \begin{equation} \label{interlacing1}
 \mu_0 \geq \nu_1 \geq \mu_2 \geq \dots \geq \mu_n \geq \nu_n \geq \mu_{n+1}=0.
\end{equation}
 
If $n=2m+1$, 
 then $\mu =(\mu_0, \cdots, \mu_{m+1})\in \Lambda^+(m+1)$
 and $\nu =(\nu_0, \cdots, \nu_{m})\in \Lambda^+(m+1)$ and 
\begin{equation}
\label{interlacing2}
  \mu_0 \geq \nu_1 \geq \mu_2 \geq \dots \geq \mu_n \geq \nu_n \geq 0.  
\end{equation}
  
 We represent the result graphically in the following theorem by representing  nontrivial symmetry breaking operators  by arrows connecting the $\theta$-stable parameters of the representations.

\begin{theorem}\label{th:branch2}
Two  representations in the standard sequences of $F^G(\mu)_\delta$ respectively $F^{G'}(\nu)_\delta$ have a nontrivial symmetry breaking operator if and only if the $\theta$-stable parameters of the representations satisfy one of the following conditions:
\\
First case: $n=2m$.  
\begin{center} 
{\rm{Case A}}
\end{center}
\begin{eqnarray*}
& (\mu_{0}, \dots , \mu_i ~|| ~\mu_{i+1} , \dots , \mu_{m+1})_{\delta}& \\
& \Downarrow &  \\
 &       (\nu_0 , \dots  ,\nu_i ~ ||~ \nu_{i+1}, \dots ,\nu_{m}  )_{\delta}
 \end{eqnarray*}
 
 \medskip
 \begin{center}
                      {\rm{or  Case B}}
  \end{center}                    
\begin{eqnarray*}
& (\mu_{0}, \dots , \mu_i ~|| ~\mu_{i+1} , \dots , \mu_{m+1})_{\delta}&  \\
& \Downarrow &  \\
&        (\nu_0 , \dots  ,\nu_{i-1} ~ ||~ \nu_i, \nu_{i+1}, \dots ,\nu_{m}  )_{\delta} &
 \end{eqnarray*}        
\medskip

Second case:  $n=2m+1$.  
\begin{center} 
{\rm{Case A}}
\end{center}
\begin{eqnarray*} 
 & (\mu_{0}, \dots , \mu_i ~|| ~\mu_{i+1} , \dots , \mu_{m+1})_{\delta} &  \\  
 & \Downarrow &    \\
 &  (\nu_0, \dots  ,\nu_i ~ ||~ \nu_{i+1}, \dots ,\nu_{m+1}  )_{\delta} &
 \end{eqnarray*}
 
\medskip

\begin{center}
{\rm{or Case B}}
\end{center}
\begin{eqnarray*}
& (\mu_{0}, \dots , \mu_i ~|| ~\mu_{i+1} , \dots , \mu_{m+1})_{\delta} &\\
& \Downarrow &  \\
 &  (\nu_0,  \dots  ,\nu_{i-1}  ~ ||~ \nu_{i}, \dots ,\nu_{m+1}  )_{\delta} &
 \end{eqnarray*}
\end{theorem}

\medskip 

\subsection{ Branching laws for representations in $\mathcal A$: \\ Third formulation}

We summarize the results as follows:

\begin{theorem}
[branching law]
\label{branching3}
Let  $\Pi$ and $\pi$  be irreducible representations in $\mathcal A$
 of $SO(n+1,1)$
 respectively $SO(n,1)$.  
\begin{enumerate}
\item[{\rm{(1)}}]
Suppose first that $n=2m$. 
Then 
\[
\operatorname{Hom}_{SO(n,1)}(\Pi|_{SO(n,1)}, \pi) \not = \{0\}
\]
 if and only if 
the enhanced $\theta$-stable parameters of the representations satisfy 
 the following conditions:
\begin{enumerate}
\item[{\rm{(a)}}] $\Pi $ and $\pi$ have the same signature $\delta$;
\item[{\rm{(b)}}]  $h_\pi \in \{ h_\Pi, h_\Pi -1\}$;
\item[{\rm{(c)}}]  
$
 \mu_0 \geq \nu_0 \geq \mu_1 \geq \dots \geq \mu_n \geq \nu_n \geq \mu_{n+1}.  
$
\end{enumerate}

\item[{\rm{(2)}}]
Suppose  that $n=2m+1$. 
Then 
\[
\operatorname{Hom}_{SO(n,1)}(\Pi|_{SO(n,1)}, \pi) \not = \{0\}
\]
 if and only if 
the enhanced $\theta$-stable parameters of the representations satisfy
 the following conditions: 
\begin{enumerate}
\item[{\rm{(a)}}] they have the same signature $\delta$;
\item[{\rm{(b)}}]  $h_\pi \in \{ h_\Pi, h_\Pi -1\}$;
\item[{\rm{(c)}}]  $
 \mu_0 \geq \nu_0 \geq \mu_1 \geq \dots \geq \mu_n \geq \nu_n \geq 0
 $.  
\end{enumerate}
\end{enumerate}
\end{theorem}

\bigskip
\section{Gross--Prasad conjectures for tempered 
 representations of 
$(SO(n+1,1), SO(n,1))$}

In this section we discuss the Gross--Prasad conjectures
 for irreducible tempered representations in $\mathcal A$ which are nontrivial on the center.  
This is a generalization of the results  for irreducible tempered representations with infinitesimal character $\rho$ which are nontrivial on the center
in \cite[Chap.~11]{sbonvec}. 
For simplicity we discuss here only the case $n=2m$.  

\medskip
Recall that for $n = 2m$,
 tempered representations $\Pi_m$ of $G= SO(n+1,1)$
 in $\mathcal A$ are irreducible unitary principal series representations  
 \cite[Prop.~15.5]{sbonvec},  
  with height $m$, signature $\delta$ and $\theta$-stable parameter
\[
   (\mu_1, \mu_2, \dots, \mu_{m}\ ||\  0 )_\delta .
\]
The tempered representations $\pi_m$ of $G'=SO(n,1)$ in $\mathcal A$ are discrete series representations with height $m$.  
Their signature is not unique and their $\theta$-stable parameter is
\[ 
   (\nu_1, \nu_2, \dots, \nu_{m}\ ||)_\delta . 
\]

We now assume that the representation $\Pi=\Pi_m$ is nontrivial on the center. These data determine  Vogan packets
 $VP({\Pi}_m)$ and $VP({\pi}_m)$ of the representations $\Pi$ of $G$,
 respectively $\pi$ of the subgroup ${G'}$.
By Theorem \ref{branching3}, 
 there is a nontrivial symmetry breaking operator 
\[ B \colon \Pi \rightarrow \pi \]
 if and only if the interlacing conditions are satisfied. 
Following exactly the steps of the algorithm by Gross and Prasad outlined
 in \cite{sbonGP}, 
 see also \cite[Chap.~11]{sbonvec}, 
we conclude that 
the Gross--Prasad conjecture predicts correctly 
 that the pair $(\Pi, \pi)$ in $(VP(\Pi_m), VP(\pi_m))$ 
 has a nontrivial symmetry breaking operator.  
 
\medskip

Together with the results about tempered  principal series representations 
 in \cite[Chap.~11, Sect.~4]{sbonvec} this completes the proof
 of the following:
\begin{theorem}
The Gross--Prasad conjectures are correct for tempered representations
 of the pair $(G,G')=(SO(n+1,1), SO(n,1))$ which are nontrivial on the center of $G$ and $G'$.
\end{theorem}

\medskip
\begin{remark}
The third formulation of the branching laws (Theorem \ref{branching3}) shows
 that branching for representations in $\mathcal A$ \lq\lq{interpolate}\rq\rq\
 between the classical branching laws of finite-dimensional representations
 and the branching laws of  Gross--Prasad for tempered representations.
\end{remark}

\bigskip

\section{Distinguished representations and periods.} 

We discuss periods of a pair of irreducible representations $\Pi$ of $G=S O(n+1,1)$
 and $\pi$ of the subgroup $G'=S O(n,1)$. 
Using the branching law (Theorem \ref{thm:branch1}), 
 we see in this section that we can prove that the representations $ A_{\mathfrak{q}}(\lambda) \in {\mathcal A}$ 
 are distinguished for some orthogonal group $H$.  

\subsection{Periods}
We recall from \cite[Thm.~5.4]{sbonvec}
 that for representations $\Pi$, $\pi$ of a real reductive Lie group $G$,
 respectively of a reductive subgroup $G'$,
 the space of symmetry breaking operators 
\[ 
     \operatorname{Hom}_{G'}(\Pi|_{G'}, \pi^{\vee}) 
\]
 and the space of $G' $-invariant continuous linear functionals
\[\operatorname{Hom}_{G'}(\Pi \boxtimes \pi,\bC)\]
 are naturally isomorphic to each other. 
Thus
 we may use symmetry breaking operators
 to construct $G'$-invariant continuous  linear functionals. 
This technique allows us to obtain
$G'$-invariant continuous linear functionals
 not only for unitary representations
 but also for nonunitary representations.  
 
 \medskip
\begin{definition}
 A nontrivial linear functional $\mathcal F$ on $\Pi \boxtimes \pi$ 
 is called a {\it period} of $\Pi \boxtimes \pi$
 if ${\mathcal{F}}$ is invariant under the diagonal $G'$-action,
 {\it{i.e.}}, 
$
{\mathcal F}  \in  \operatorname{Hom}_{G'}(\Pi \boxtimes \pi,\bC).  
$  
\end{definition}
 
We say that vector
 $\Phi \otimes \phi \in \Pi \boxtimes \pi$ is a {\it{test vector}} for the period $\mathcal F$ 
 if $\Phi \otimes \phi$ is not in the kernel of ${\mathcal F}$. 
If the period is nontrivial on a  test vector  
 $\Phi \otimes \phi$,
 we refer to its image as the {\it value of the period} on $\Phi \otimes \phi$. 
 
\medskip
\begin{remark}
If $(\Pi, \pi)$ is a pair of discrete series representations  for the symmetric pairs $(G_1(\mathbb R),G_2(\mathbb R))$
 we may consider a realization of $\Pi \boxtimes \pi$
 in $L^2( G_1(\mathbb R) \times G_2(\mathbb R)) $.
The integral 
\[  \int_{G_2(\mathbb R )}
\Phi (h) \phi (h) dh
\]
converges for some  smooth vectors $(\Phi,\phi ) \in \Pi \boxtimes \pi$ in the minimal $K$-types and so if it is nonzero 
it defines a  period integral for the discrete series representations $\Pi \boxtimes \pi$ \cite{Vargas}.
If the representation $\Pi \boxtimes \pi$ is not tempered
 the integral usually does not converge,
 but nevertheless we can consider periods
 via symmetry breaking operators.  
\end{remark}

The next theorem describes 
 for the pair $(G,G')= (SO(n+1,1), SO(n,1))$ the $\theta$-stable parameters of the representations in $\mathcal A$ which have a nontrivial period \linebreak
  ${\operatorname{Hom}}_{G'}(\Pi \boxtimes \pi, {\mathbb{C}}) .$
Recall that the $\theta$-stable parameter of a representation
 $\Pi \in {\mathcal A}$ of height $i$ is of the form 
\begin{enumerate}
\item $ (\mu_{1}, \dots , \mu_i ~|| ~\mu_{i+1} , \dots  ,0)_{\delta} $ if $G=SO(2m+1,1)$, 
\item
$ (\mu_{1}, \dots , \mu_i ~|| ~\mu_{i+1} , \dots  ,\mu_m)_\delta$ if $G'=SO(2m,1)$.  
\end{enumerate}

\medskip
\begin{theorem}
\label{thm:branch4}
Suppose that  $\Pi $ and $\pi$  are representations  of $G$,
 respectively $G'$ in $\mathcal A$ 
 of  height $i$ respectively j. 
The following statements are equivalent:
\begin{enumerate}
\item[{\rm{(i)}}] The representation $\Pi \boxtimes \pi$ has a nontrivial $G'$-period;
\item[{\rm{(ii)}}] $\Pi $ and $\pi$ have the same signature, 
 $j=i$ or $i-1$ and their $\theta$-stable parameters satisfy the interlacing conditions of the branching result in Theorem \ref{branching3}.  
\end{enumerate}
\end{theorem}

\medskip

In \cite[Props.~10.12 and 10.30]{sbonvec}
 we proved furthermore for $(O(n+1,1), O(n,1))$.

\begin{theorem}
\label{thm:testvector}
If $\Pi$ and $\pi$ are representations of O(n+1,1) respectively O(n,1) with trivial infinitesimal character $\rho$
such that $\Pi \otimes \pi$ has a nontrivial $G'$-period.  
Then
 there is a test vector for a nonzero period in the minimal $K$-type 
$\Pi \boxtimes \pi$.  
\end{theorem}

\medskip
\begin{remark}  We expect that Theorem \ref{thm:testvector} also holds for unitary representations in $\mathcal A$.
\end{remark}

\medskip
\begin{remark}
Similar results for cohomologically induced  representations
 of other pairs $(G,G')$ of reductive groups were obtain by B.~Sun \cite{sun2}.  
\end{remark}

\medskip
\subsection{Distinguished representations}
Let $G$ be a reductive group, 
 and $H$ a reductive subgroup.  
We regard $H$ as a subgroup of the direct product group $G \times H$
 via the diagonal embedding $H \hookrightarrow G \times H$.  

\begin{definition}
Let $\psi$ be a one-dimensional representation of $H$. 
We say an admissible smooth representation  $\Pi$ of $G$ is
$(H,\psi)$-{\it{distinguished}}
 if 
\[ \operatorname{Hom}_{H}(\Pi \boxtimes \psi^{\vee}, {\mathbb{C}})
 \simeq \operatorname{Hom}_{H}(\Pi|_H, \psi) \not = \{0\}.  
\]
If the character $\psi $ is trivial,
 we say $\Pi$ is $H$-{\it{distinguished}}.  
\end{definition}

\medskip

Repeated application of Theorem \ref{branching3} proves

\begin{theorem}
Suppose that $\Pi =A_\bq(\lambda) \in {\mathcal A}$ is a representation of $S O(n+1,1)$
 of height h. 
Then $\Pi $ is $S O(n+1-h,1)$-distinguished.
\end{theorem}

\medskip
\begin{remark}
For a different proof and perspective of this theorem see \cite{K1}.
\end{remark}

\medskip 
Since the representations $\Pi_{i,+}$ have height $i$
 (see Example \ref{ex:reprho}), 
 this generalizes the following theorem proved 
in \cite[Thm.~12.4]{sbonvec}. 

\begin{theorem}
Let $0 \le i \le \frac n 2$.  
Then the representations $\Pi_{i, \delta}$ $(\delta \in \{+, -\})$
 of $G=SO(n+1,1)$ are  $SO(n+1-i,1)$-distinguished. 
\end{theorem}

\medskip
In the remainder of the section 
 we recall a formula for the period of representations with trivial infinitesimal character $\rho $ of the pair
\[
   (G,H)=(SO(n+1,1), SO(m+1,1))
\quad
\text{for  $m \le n$}.
\]
We use here the notation of Example \ref{ex:reprho} in Section \ref{subsec:II.2.2}.

The period can be computed  by applying  the composition 
 of the regular symmetry breaking operators
 that we constructed in \cite[Chap.~12, Sect.~1]{sbonvec}
 with respect to the chain
 of subgroups
\begin{equation}
\label{eqn:seqOn}
  G=SO(n+1,1) \supset SO(n,1) \supset SO(n-1,1) \supset \cdots \supset SO(m+1,1)=H,
\end{equation}
  to test vectors. 
We write simply $\Pi_i$ for $\Pi_{i,+}$.  
We recall from \cite[Prop.~14.44]{sbonvec}
 that $\Pi_i \equiv \Pi_{i,+}$ has a minimal $K$-type
 $K_{min}(i,+)={\bigwedge^i(\mathbb{C}}^{n+1}) \boxtimes {\bf{1}}$.

Let $v\in \bigwedge^i({\mathbb{C}}^{n+1})$ be the image of $1 \in {\mathbb{C}}$
 via the following successive inclusions:
\[
   {\Exterior}^{i}({\mathbb{C}}^{n+1}) \supset 
   {\Exterior}^{i-1}({\mathbb{C}}^{n}) \supset 
   \cdots \supset
   {\Exterior}^{i-\ell}({\mathbb{C}}^{n+1-\ell}) \supset 
   \cdots \supset 
   {\Exterior}^{0}({\mathbb{C}}^{n+1-i}) \simeq {\mathbb{C}} \ni 1, 
\]
 and we regard $v$ as an element
 of the minimal $K$-type $K_{min}(i,+)$
 of $\Pi_i$.

\begin{theorem}
[{\cite[Thm.~12.5]{sbonvec}}]
\label{thm:period2}
Let $\Pi_i$ be the irreducible representation of $G=SO(n+1,1)$,
 with infinitesimal character $\rho$, 
 height $i$ and signature $+$. 
Let 
 $v$ be the normalized element
 of its minimal $K$-type as above.  
For $0 \le i \le n$, 
 the value $F(v)$ of the $SO(n+1-i,1)$-period $F$ 
 on $v \in \Pi_i$ is
\[
  \frac{\pi^{\frac 1 4 i(2n-i-1)}}{((n-i)!)^{i-1}}
  \times
\begin{cases}  
  \frac{1}{(n-2i)!}\quad &\text{if $2i < n+1$}, 
\\
  (-1)^{n+1} (2i-n-1)!\quad & \text{if $2i \ge n+1$}.  
\end{cases}
\]
\end{theorem}

\medskip

\section{Bilinear forms on $(\bg,K)$-cohomologies induced by symmetry breaking operators}
\label{sec:6}
Consider now the induced map by a symmetry breaking operator
\[ T\colon \Pi \rightarrow \pi\]  on  $(\bg,K)$-cohomologies  of a pair of representations $\Pi$ and $\pi$.  In what follows,
 by abuse of notation,
 we denote an admissible  smooth representation and its underlying $({\mathfrak{g}},K)$-module by the same letter 
 when we discuss their $({\mathfrak{g}},K)$-cohomologies.

\bigskip
Recall that  a theorem of Vogan--Zuckerman \cite{VZ} states that every irreducible {\it unitary} representation $\Pi $ of $SO(n+1,1)$ with 
\[H^*(\bg,K;\Pi\otimes V) \not = \{0\}\] 
for a finite-dimensional representation $V$ is of the form $\Pi = A_{\mathfrak{q}}(\lambda)$. If we assume that 
\[H^*(\bg,K;\Pi) \not = \{0\}\]
  then $\Pi $ is isomorphic  to a unitary irreducible representation with infinitesimal character $\rho$ 
{\it{i.e.}},  
 it is of the form $A_{\mathfrak{q}}$. 
See \cite[Chap.~14, Sect.~9.4]{sbonvec} for $O(n+1,1)$.
 
Note also that an irreducible representation $\Pi$ with 
\[H^*(\bg,K;\Pi\otimes V) \not = \{0\}\]
 for some finite-dimensional representation $V$
 is not always unitarizable.  

\medskip
Suppose that  $\Pi $ is a principal series representation of a connected reductive Lie group with nonsingular integral infinitesimal character. 
If the $(\bg,K$)-cohomology of $\Pi \otimes V$ is nonzero the highest weight of $V$ satisfies the conditions in \cite[Chap.~III, Thm.~3.3]{BW}.  
For representations of $O(n+1,1)$  the situation is  more complicated  and the finite-dimensional representation $V$ is also described
 in \cite[Chap.~16, Sect.~4]{sbonvec}. 
Using the results in \cite[Chap.~15]{sbonvec} about the restriction of representations of $O(n+1,1)$ to $SO(n+1, 1)$ we obtain a formula for the representation $V$ for $SO(n+1,1)$.

\bigskip
Let $(G,G')=(SO(n+1,1),SO(n,1))$, 
 and $\Pi$, $\pi$ be representations 
 of $G$ and $G'$, 
 respectively with 
\[H^*(\bg,K;\Pi\otimes V) \not = \{0\}\]
  and 
\[H^*(\bg',K'; \pi\otimes V') \not = \{0\}, \]
where $V$ and $V'$ are irreducible finite-dimensional representations
 of $G$ and $G'$, 
respectively.  
Suppose in addition that 
\begin{enumerate}
\item $\mbox{Hom}_{G'}(V|_{G'},V') \not = \{0\}$;
\item $\Pi$ and $\pi$ have the same height $i$;
\item  $\Pi$ and $\pi$ have the same signature $\delta$.  
\end{enumerate}
A nontrivial symmetry breaking operator 
 $T\colon \Pi \otimes V \to \pi \otimes V'$  induces a canonical homomorphism 
 \begin{equation}
T^{\ast} \colon H^j({\mathfrak{g}}, K; \Pi \otimes V) \to H^j({\mathfrak{g}}', K'; \pi \otimes V')
\end{equation}
 and a bilinear form 
\[
   B_T \colon 
   H^j({\mathfrak{g}}, K; \Pi \otimes V )
   \times 
   H^{n-j}({\mathfrak{g}}', K'; (\pi \otimes V')^{\vee}\otimes \chi_{(-1)^{n}})
   \to {\mathbb{C}}
\quad
\text{for all $j$.}
\]
where $(\pi \otimes V')^{\vee}$ denotes the contragredient representation
 of $\pi \otimes V'$. 

\bigskip

The formulas in \cite[Chap.~16, Sect.~3]{sbonvec}
 and \cite[Chap.~III, Thm.~3.3]{BW} imply
\begin{theorem}
\label{thm:Bnonzero}
 Suppose that $\Pi$ and $\pi$ are principal series representations and that 
\begin{enumerate}
\item $H^*(\bg,K, \Pi \otimes V) \not = \{0\}$
 and $H^*(\bg',K', \pi \otimes V') \not = \{0\}$;
\item $\operatorname{Hom}_{G'}(V|_{G'},V') \not = \{0\}$;
 \item $\Pi$ and $\pi$ have the same height $i$; 
 \item  $\Pi$ and $\pi$ have the same signature $\delta$;
 \item there exists a nontrivial symmetry breaking operator $T\colon \Pi \to \pi $.  
 \end{enumerate}
The symmetry breaking operator $T$ induces 
a nontrivial homomorphism
\[
   T^i \colon
   H^i({\mathfrak{g}}, K; \Pi \otimes V )
   \rightarrow
   H^{i}({\mathfrak{g}}', K'; \pi \otimes V')
\] 
and hence a nontrivial bilinear form

\[
   B_T \colon H^i({\mathfrak{g}}, K; \Pi \otimes V ) \times 
              H^{n-i}({\mathfrak{g}}', K'; (\pi \otimes V')^\vee\chi_{(-1)^{n}} ) \rightarrow \bC. \]

\end{theorem}

\medskip
The following theorem provides a criterium for the nonvanishing of this bilinear form on the $(\bg,K)$-cohomology of representations with trivial infinitesimal character.

\begin{theorem}
[{\cite[Thm.~12.11]{sbonvec}}]
\label{thm:171556}
Let $T \colon X \to Y$ be a $({\mathfrak{g}}',K')$-homomorphism, 
 where $X$ is a $({\mathfrak{g}},K)$-module $A_{\mathfrak{q}}$
 and $Y$ is a $({\mathfrak{g}}',K')$-module $A_{\mathfrak{q}'}$.  
Let $U$ be the representation space
 of the minimal $K$-type $\mu$ in $X$, 
 and $U'$ that of the minimal $K'$-type $\mu'$ in $Y$.  
We define a $K'$-homomorphism by 
\begin{equation}
   \varphi_T := {\operatorname{pr}} \circ T|_{U}
   \colon U \to U'.  
\end{equation}
\begin{enumerate}
\item[{\rm{(1)}}]
If $\varphi_T$ is zero, 
 then the homomorphisms 
$
  T_{\ast} \colon 
  H^j({\mathfrak{g}},K;X) \to H^j({\mathfrak{g}}',K';Y)
$
 and the bilinear form $B_T$
 vanish for all degrees $j \in {\mathbb{N}}$.   
\item[{\rm{(2)}}]
If $\varphi_T$ is ${\mathfrak{p}}$-nonvanishing
 at degree $j$, 
 then $T_{\ast}$ and the bilinear forms $B_T$ are nonzero
  for this degree $j$.  
\end{enumerate}
\end{theorem}

This theorem together with our results \cite{sbonvec} implies

\begin{theorem}
[cf.~ {\cite[Thm.~12.13]{sbonvec}}]
\label{thm:bilinear}
Let $(G,G')=(SO(n+1,1),SO(n,1))$, 
 $0 \le i \le n$, 
 and $\delta \in \{ +,- \}$.  
Let $T$ be the symmetry breaking operator 
 $\Pi_{i,\delta} \to \pi_{i,\delta}$
 given  in Theorem \ref{thm:branch1}  
\begin{enumerate}
\item[{\rm{(1)}}]
$T$ induces bilinear forms
\[
   B_T \colon 
   H^j({\mathfrak{g}}, K; \Pi_{i,\delta}) 
   \times 
   H^{n-j}({\mathfrak{g}}', K'; \pi_{n-i,(-1)^n \delta})
   \to {\mathbb{C}}
\quad
\text{for all $j$.}
\]
\item[{\rm{(2)}}]
The bilinear form $B_T$ is nonzero
 if and only if $j=i$ and $\delta=(-1)^i$.  
\end{enumerate}
\end{theorem}

\medskip

\begin{remark}
A theorem similar  to Theorem \ref{thm:bilinear}
 was proved 
 by B.~Sun \cite{sun}
 for the $(\mathfrak{g},K)$-cohomology with nontrivial coefficients
 of irreducible tempered representations
 of the  pair \[(GL(n,\mathbb R), GL(n-1,\mathbb R)).\]
\end{remark}

\smallskip
{\bf{Acknowledgements}}\enspace
The first author was partially supported
 by Grant-in-Aid for Scientific Research (A) (JP18H03669), 
Japan Society for the Promotion of Science.
The second author was  partially supported by NSF grant DMS-1500644.

\end{document}